\documentclass[11pt]{article}
\usepackage{amsmath,amssymb,psfig,epsfig,latexsym,graphicx,here}

\newcommand{\Snj}{\scriptstyle}
\newcommand{\beql}[1]{\begin{equation}\label{#1}}
\newcommand{\eeq}{\end{equation}}
\newcommand{\eqn}[1]{(\ref{#1})}
\newcommand{\F}{\ensuremath{\frac}}

\begin{document}
\begin{center}
{\large\bf Seven Staggering Sequences} \\
\vspace*{+.1in}
{\em N. J. A. Sloane} \smallskip \\
Algorithms and Optimization Department  \\
AT\&T Shannon Lab \\
Florham Park, NJ 07932--0971 \medskip \\
Email address: {\tt njas@research.att.com}
\bigskip \\
April 3, 2006 \bigskip
\end{center}

\setlength{\baselineskip}{1.0\baselineskip}

\paragraph{0. Introduction}
When the {\em Handbook of Integer Sequences} came out in 1973,
Philip Morrison gave it an enthusiastic review in the {\em Scientific
American} and Martin Gardner 
was kind enough to say in his {\em Mathematical Games} column
for July 1974 that ``every recreational mathematician should buy
a copy forthwith.''
That book contained 2372 sequences.
Today the {\em On-Line Encyclopedia of Integer Sequences}
(or OEIS)
\cite{OEIS} contains 117000 sequences.
The following are seven that I find especially interesting.
Many of them quite literally stagger.
The sequences will be labeled with their numbers (such
as A064413) in the OEIS.
Much more information about them can be found there and
in the references cited.

\paragraph{1. The EKG sequence}
(A064413, due to Jonathan Ayres).
The first three sequences are defined by unusual
recurrence rules.
The first begins with $a(1)=1, a(2)=2$, and the rule
for extending it is that the next term, $a(n+1)$, is taken to be
the smallest positive number not already
in the sequence which has a nontrivial common factor
with the previous term $a(n)$.
Since $a(2)=2$, $a(3)$ must be even, and is therefore $4$;
$a(4)$ must have a factor in common with $4$, that is, must also
be even, and so $a(4)=6$.
The smallest number not already in the sequence that has a common factor
with $6$ is $3$, so $a(5)=3$, and so on.
The first 18 terms are
$$
1, 2, 4, 6, 3, 9, 12, 8, 10, 5, 15, 18, 14, 7, 21, 24, 16, 20, \ldots \, .
$$
It is clear that if a prime $p$ appears in the sequence,
$2p$ will be the term either immediately before or after it.
Jeffrey Lagarias, Eric Rains and I studied this
sequence in \cite{EKG}.
One of the things that we observed was that in fact
every odd prime $p$ was always {\em preceded} by
$2p$, and always {\em followed} by $3p$. This is certainly true for the
first $10000000$ terms, but we were unable to prove it in general.

\begin{figure}[htb]
\centerline{\includegraphics[angle=90, width=4in]{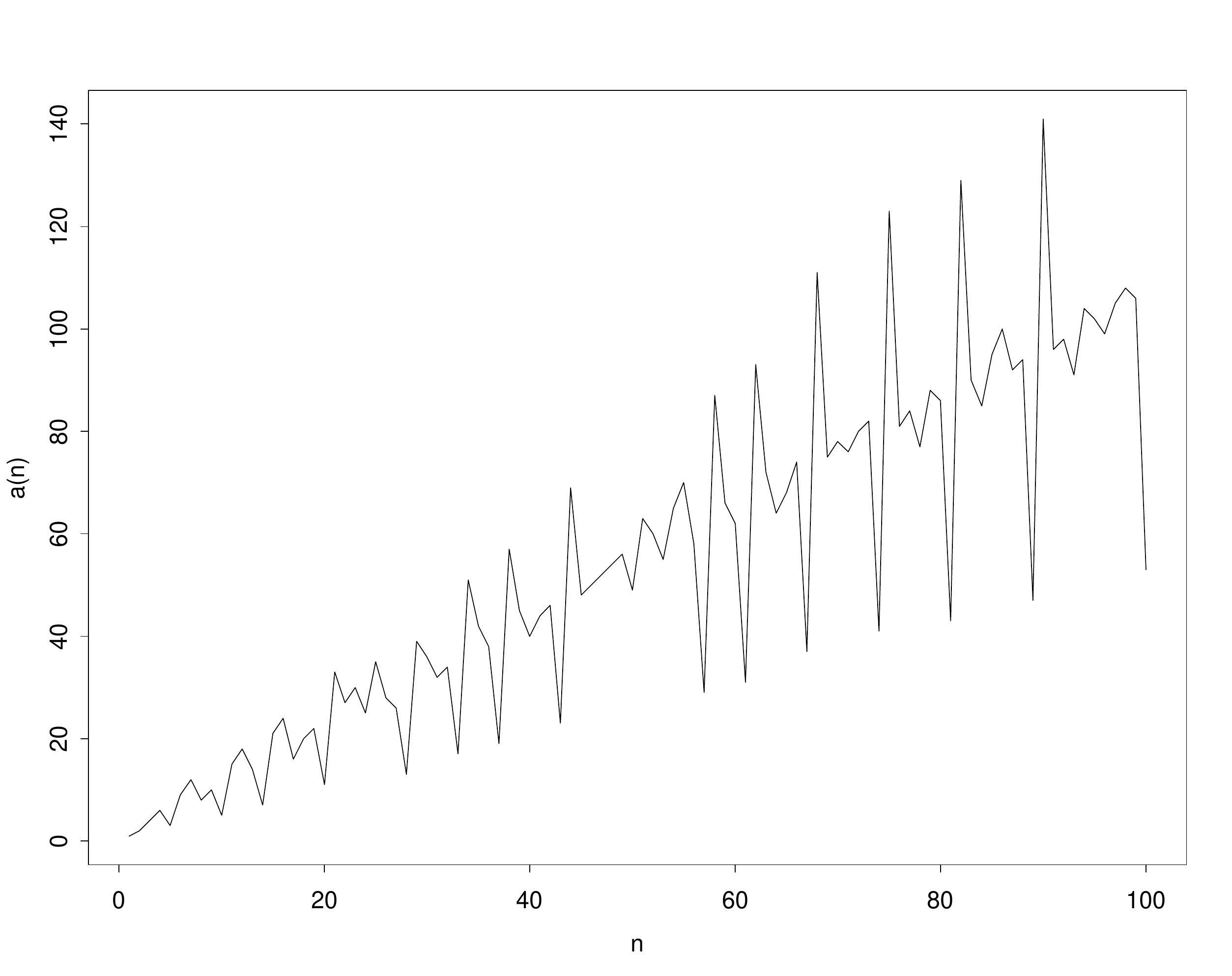}}

\caption{The first 100 terms of the EKG sequence,
with successive points joined by lines.}
\label{F1}
\end{figure}

\begin{figure}[H]
\centerline{\includegraphics[angle=90, width=4in]{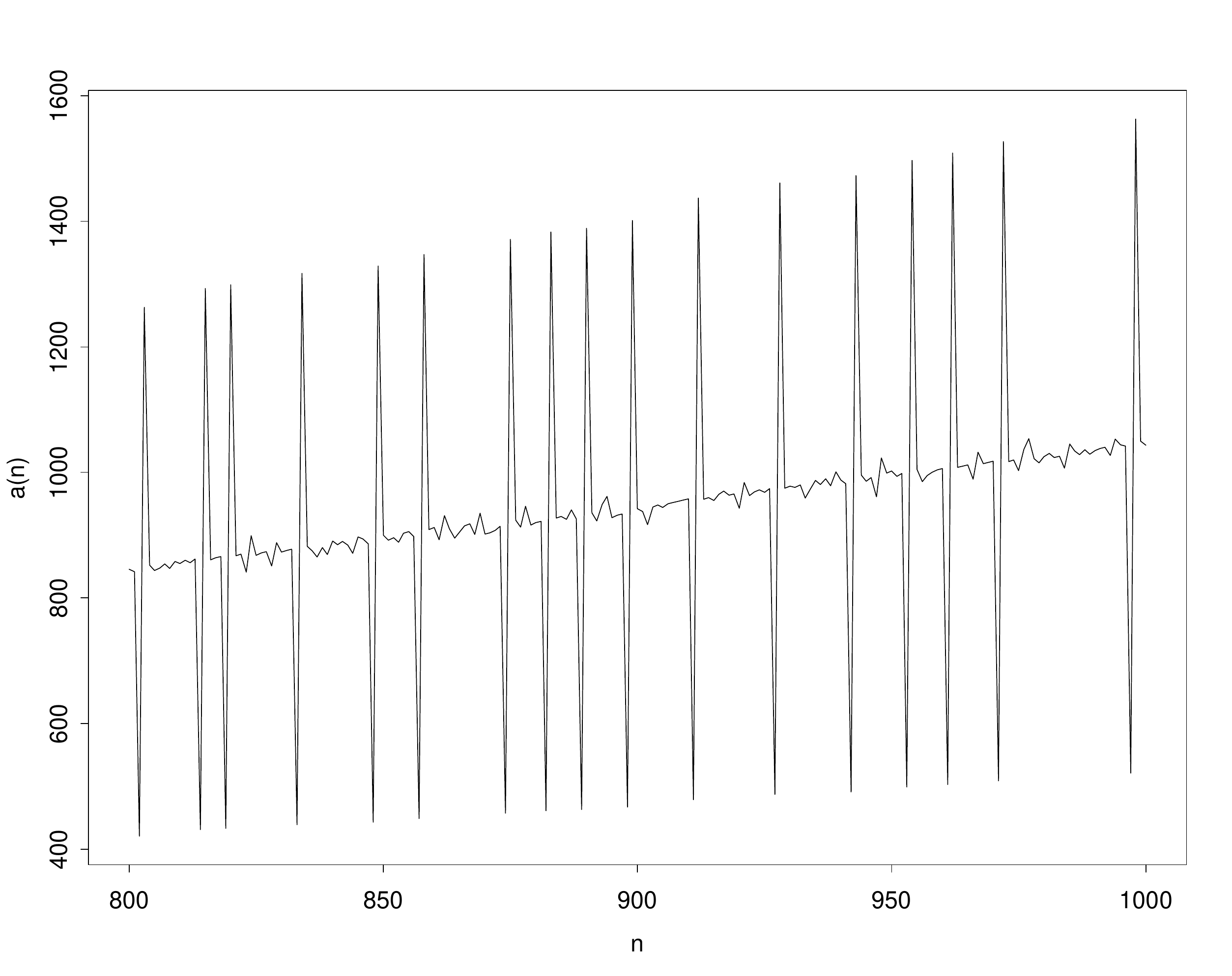}}

\caption{Terms 800 to 1000 of the EKG sequence.}
\label{F2}
\end{figure}

We called this the EKG sequence, since it looks
like an electrocardiogram when plotted
(Figs. \ref{F1}, \ref{F2}).

There is an elegant three-step proof that every positive number must
eventually  appear in the sequence.
(i) If infinitely many multiples of some prime $p$ appear
in the sequence, then every multiple of $p$ must appear.
(For if not, let $kp$ be the smallest missing multiple of $p$.
Every number below $kp$ either appears or it doesn't,
but once we get to a multiple of $p$ beyond all those terms,
the next term must be $kp$, which is a contradiction.)
(ii) If every multiple of a prime $p$ appears, then every number appears.
(The proof is similar.)
(iii) Every number appears. (For if there are only finitely many different
primes among the prime factors of all the terms, then some prime must
appear in infinitely many terms, and the result follows from (i) and (ii).
On the other hand, if infinitely many different primes $p$ appear,
then there are infinitely many numbers $2p$, as noted above,
so $2$ appears infinitely often, and again the result
follows from (i) and (ii).)

\begin{figure}[htb]
\centerline{\includegraphics[angle=0, width=4in]{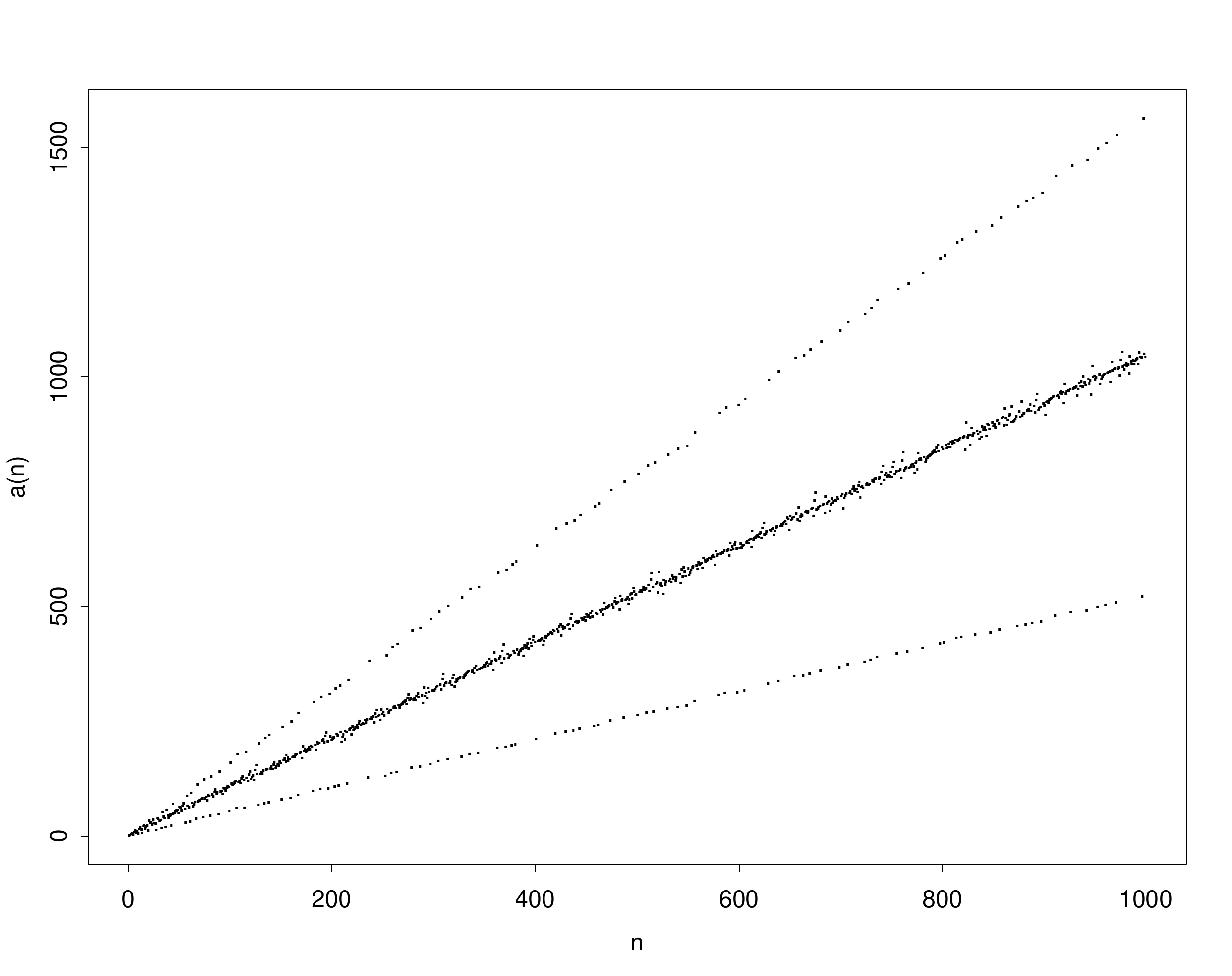}}

\caption{The first 1000 terms of the EKG sequence, successive points not joined.
They lie roughly on three almost-straight lines.}
\label{F3}
\end{figure}

Although the initial terms of the sequence stagger around,
when we look at the big picture we find that the points lie
very close to three almost-straight lines
(Fig. \ref{F3}).
This is somewhat similar to the behavior of the prime numbers,
which are initially erratic, but lie close to a smooth curve
(since the  $n^{\rm th}$ prime is roughly $n \log n$)
when we look at the big picture---see Don Zagier's
lecture on ``The first 50 million prime numbers'' \cite{Zagier77}.

In fact, we have a precise conjecture about the three lines
on which the points lie.
We believe---but were unable to prove---that
almost all
$a(n)$ satisfy the asymptotic formula
$a(n) \sim n (1+ 1/(3 \log n)) $ (the central line in Fig. \ref{F3}),
and that the exceptional values $a(n)=p$ and $a(n)= 3p$, for
$p$ a prime, produce the points on the lower and upper lines.
We {\em were} able to show that 
the sequence has essentially linear growth
(there are constants $c_1$ and $c_2$ 
such that $c_1 n <a(n) < c_2 n$ for all $n$),
but the proof of even this relatively weak result was quite difficult.
It would be nice to know more about this sequence!

\paragraph{2. Gijswijt's sequence} (A090822, invented
by Dion Gijswijt when he was a graduate student
at the University of Amsterdam,
and analyzed by him,
Fokko van de Bult,
John Linderman,
Allan Wilks and myself \cite{GIJ}.)
We begin with $b(1)=1$.
The rule for computing the next term, $b(n+1)$,
is again rather unusual. We write the sequence of numbers
we have seen so far,
$$
b(1), b(2), \ldots, b(n) \,,
$$
in the form of an initial string $X$, say (which can be
the empty string $\emptyset$),
followed by as many repetitions as possible
of some nonempty string $Y$. That is, we write
\beql{Eq1}
b(1), b(2), \ldots, b(n) ~=~ XY^k, 
\mbox{~where~} k \mbox{~is~as~large~as~possible} \,. 
\eeq
Then $b(n+1)$ is $k$.

Some examples will make this clear.
The sequence begins:
$$
1,1,2,1,1,2,2,2,3,1,1,2,1,1,2,2,2,3,2,1,1,2, \ldots \,.
$$
After the first six terms we have
$$
b(1), b(2), \ldots, b(6) ~=~ 1,1,2,1,1,2,
$$
so we can take $X$ to be empty, $Y$ to be $1,1,2$ and $k=2$,
so $b(1), b(2), \ldots, b(6) = Y^2$. This is the largest
$k$ we can achieve here, so $b(7) = 2$.
Now we have
$$
b(1), b(2), \ldots, b(7) ~=~ 1,1,2,1,1,2,2,
$$
and we can take $X=1,1,2,1,1$, $Y=2$, $k=2$, getting $b(8)=2$.
Next,
$$
b(1), b(2), \ldots, b(8) ~=~ 1,1,2,1,1,2,2,2,
$$
and we can take $X=1,1,2,1,1$, $Y=2$, $k=3$, getting $b(9)=3$,
the first time a $3$ appears. And so on.

The first time a $4$ appears is at $b(220)$. We computed
several million terms  without finding a $5$,
and for a while we wondered if perhaps no term greater than $4$
was ever going to appear.
However, we were able to show that a $5$ does
eventually appear, although the universe would grow
cold before a direct search would find it.
The first $5$ appears at about term 
$$
10^{10^{\Snj 23}} \, .
$$
The sequence is in fact unbounded, and the first time that
a number $m~(= 5,6,7,\ldots)$ appears seems to be at about term number
$$
2^{2^{\Snj 3^{\Snj 4^{\cdot^{\cdot^{\cdot^{\Snj {m-1}}}}}}}} \,,
$$
a tower of height $m-1$.

There are of course several well-known sequences
which have an even slower growth-rate than this one
(the inverse Ackermann function \cite{Ack28},
the Davenport-Schinzel sequences \cite{ShAg95},
or the inverse to Harvey Friedman's sequence \cite{Fri01}, for example). 
Nevertheless, I think the combination of slow growth and
an unusual definition make Gijswijt's sequence remarkable.
It also has an interesting recursive structure,
which is the key to its analysis.
There is only room here to give a hint of this.

The starting point is the observation that the
sequence---let's call it $A^{(1)}$---can be built up recursively from 
``blocks'' that are always doubled and are followed by ``glue'' strings.
The first block is $B_1 = 1$, the first glue string is 
$S_1 = 2$, and the sequence begins with
$$
B_1 B_1 S_1 ~=~ 1, 1, 2 \,,
$$
which is the second block, $B_2$.
The second glue string is 
$S_2 = 2, 2, 3$, and the sequence also begins with
$$
B_2 B_2 S_2 ~=~ 1, 1, 2, 1, 1, 2, 2, 2, 3 \,,
$$
which is the third block, $B_3$.
This continues: for all $m$, the sequence
begins with $B_{m+1} = B_m B_m S_m$,
where $S_m$ contains no $1$'s and is terminated
by the first $1$ which follows $B_m B_m$.
Now something remarkable happens: if we concatenate 
all the glue strings $S_1, S_2, S_3, \ldots$,
we get a new sequence, $A^{(2)}$ say:
$$
2, 2, 2, 3, 2, 2, 2, 3, 2, 2, 2, 3, 3, 2, 2, 2, 3, 2, 2, 2, 3, 2, 2, 2, 3, 3, 2, 2, 2, 3, 2, 2, 2, 3, 2, 2, 2, 3, 3, 3, 3, \ldots 
$$ 
(A091787), which turns out to generated by the same rule, \eqn{Eq1}, as the original sequence,
except that the next term is now the maximum of 
$k$ and $2$. If $k=1$ is the best we
can achieve, we promote it to $2$.
We call $A^{(2)}$ the {\em second-order} sequence.
This has a similar recursive structure to the original sequence, only
now it is built up from blocks which are repeated three times
and followed by second-order glue strings which contain no $1$'s or $2$'s.
If we concatenate the second-order glue strings we get
the {\em third-order} sequence $A^{(3)}$, which 
is built up from blocks which are repeated four times
and followed by third-order glue strings which contains no $1$'s, $2$'s
or $3$'s. And so on.

Now we observe that arbitrarily long initial segments of
the second-oder sequence  $A^{(2)}$ appear as
subsequences of the original sequence $A^{(1)}$,
arbitrarily long initial segments of
the third-oder sequence  $A^{(3)}$ appear as
subsequences of the second-order sequence $A^{(2)}$, $\ldots$.
But the $m^{\rm th}$-order sequence $A^{(m)}$ begins with $m$.
So the original sequence contains every positive number!

Of course all this requires proof, and the reader is referred
to \cite{GIJ} for further information about this fascinating sequence.


We observed experimentally that in
variations of Gijswijt's sequence with initial 
conditions consisting of any finite string of $2$'s and $3$'s, a $1$ always 
eventually appeared in the sequence, but were unable to prove that this 
would always be the case.
We called this the ``Finiteness Conjecture'':
start with any finite initial string of numbers,
and extend it by the ``next term is $k$'' rule \eqn{Eq1}.
Then eventually one must see a $1$.
If we had a direct proof of this, it would simplify the 
analysis of the original sequence.  Can some reader find a proof?

\paragraph{3. Numerical analogs of Aronson's sequence.}
Aronson's sentence is a classic self-referential assertion:
``$t$ is the first, fourth, eleventh, sixteenth, $\ldots$ letter in this sentence''
(\cite{Aro85}, \cite{Hof85})
and produces the sequence
$1, 4, 11, 16, 24, 29, 33, 35, 39, 45, \ldots$ (A005224).
It suffers from the drawback that later terms are ill-defined,
because of the ambiguity in the English names for numbers---some
people say ``one hundred and one,'' others ``one hundred one,'' etc.

Another well-known self-referential sequence is Golomb's sequence,
which is defined by the property that the $n^{\rm th}$
term is the number of times $n$ appears in the sequence
$$
1,2,2,3,3,4,4,4,5,5,5,6,6,6,6,7,7,7,7,8, \ldots 
$$
(A001462).  There is a simple formula for the  $n^{\rm th}$ term: it is
the nearest integer to (and approaches)
$$
\phi^{2- \phi} \, n^{\phi -1} \,,
$$
where $\phi = (1+ \sqrt{5})/2$ is
the golden ratio (\cite{Gol66}, \cite[Section E25]{UPNT}).

In \cite{CSV03}, Benoit Cloitre, Matthew Vandermast and I
studied some new kinds of self-referential sequences, one of which is:
$c(n+1)$ is the smallest positive number $> c(n)$ consistent with the condition
``$n$ is a member of the sequence if and only if $c(n)$ is odd.''

What is $c(1)$? Well, $1$ is the smallest positive number 
consistent with the conditions, so $c(1)$ must be $1$.
What about $c(2)$? It must be at least $2$, and it can't be $2$,
for then $2$ would be in the sequence, but $c(2)$ would be even.
Nor can it be $3$, for then $2$ would be missing (the sequence
increases) whereas $c(2)$ would be odd.
But $c(2)$ {\em could} be $4$, and therefore {\em must} be $4$.
So $c(3)$ must be even and $>4$, and $c(3) =6$ works.
Now 4 is in the sequence, so $c(4)$ must be odd, and $c(4) =7$ works.
Continuing in this way we find that the first few terms are as follows
(this is A079000):
$$
\begin{array}{rrrrrrrrrrrrrr}
n: & 1 & 2 & 3 & 4 & 5 & 6 & 7 & 8 & 9 & 10 & 11 & 12 & \cdots \\
c(n): & 1 & 4 & 6 & 7 & 8 & 9 & 11 & 13 & 15 & 16 & 17 & 18 & \cdots
\end{array}
$$

Once we are past $c(2)$ there are no further complications,
$c(n-1)$ is greater than $n$, and we {\em can},
and therefore {\em must}, take
$$
c(n) = c(n-1) + \epsilon ~,
$$
where $\epsilon$ is 1 or 2 and is given by:
$$
\begin{array}{ccc}
~ & \mbox{$c(n-1)$ even} & \mbox{$c(n-1)$ odd} \\ [+.1in]
\mbox{$n$ in sequence} & 1 & 2 \\
\mbox{$n$ not in sequence} & 2 & 1
\end{array}
$$
The gap between successive terms for $n \ge 3$ is either 1 or 2.

The analogy with Aronson's sequence is clear.
Just as Aronson's sentence indicates exactly which of its terms are t's,
$\{c(n)\}$ indicates exactly which of its terms are odd.

It is easy to show that all odd numbers $\ge 7$ occur in the sequence.
For suppose some number $2t+1$ were missing.
Therefore $c(i) = 2t$, $c(i+1) = 2t+2$ for some $i \ge 3$.  From
the definition, this means $i$ and $i+1$ are missing,
implying a gap of at least 3, a contradiction to what we just observed.

Table \ref{T1} shows the first 72 terms, with the even numbers underlined.
\begin{table}[htb]
$$
\begin{array}{rrrrrrrrrrrr}
n: & 1 & \multicolumn{1}{r|}{2} & 3 & 4 & 5 & 6 & 7 & \multicolumn{1}{r|}{8} & 9 & 10 \\
c(n): & 1 & \multicolumn{1}{r|}{\underline{4}} & \underline{6} & 7 & \underline{8} & 9 & 11 & \multicolumn{1}{r|}{13} & 15 & \underline{16} \\
~ & ~ \\
n: & 11 & 12 & 13 & 14 & 15 & 16 & 17 & 18 & 19 & \multicolumn{1}{r|}{20} \\
c(n): & 17 & \underline{18} & 19 & \underline{20} & 21 & 23 & 25 & 27 & 29 & \multicolumn{1}{r|}{31} \\
~ & ~ \\
n: & 21 & 22 & 23 & 24 & 25 & 26 & 27 & 28 & 29 & 30 \\
c(n): & 33 & \underline{34} & 35 & \underline{36} & 37 & \underline{38} & 39 & \underline{40} & 41 & \underline{42} \\
~ & ~ \\
n: & 31 & 32 & 33 & 34 & 35 & 36 & 37 & 38 & 39 & 40 \\
c(n): & 43 & \underline{44} & 45 & 47 & 49 & 51 & 53 & 55 & 57 & 59 \\
~ & ~ \\
n: & 41 & 42 & 43 & \multicolumn{1}{r|}{44} & 45 & 46 & 47 & 48 & 49 & 50 \\
c(n): & 61 & 63 & 65 & \multicolumn{1}{r|}{67} & 69 & \underline{70} & 71 & \underline{72} & 73 & \underline{74} \\
~ & ~ \\
n: & 51 & 52 & 53 & 54 & 55 & 56 & 57 & 58 & 59 & 60 \\
c(n): & 75 & \underline{76} & 77 & \underline{78} & 79 & \underline{80} & 81 & \underline{82} & 83 & \underline{84} \\
~ & ~ \\
n: & 61 & 62 & 63 & 64 &  65 & 66 & 67 & 68 & 69 & 70 \\
c(n): & 85 & \underline{86} & 87 & \underline{88} & 89 & \underline{90} & 91 & \underline{92} & 93 & 95 \\
~ & ~ \\
n: & 71 & 72 & \multicolumn{2}{l}{\cdots} \\
c(n): & 97 & 99 & \multicolumn{2}{l}{\cdots} \\
\end{array}
$$
\caption{The first 72 terms of the sequence ``$n$ is in the sequence if and only if $c(n)$ is odd.''}
\label{T1}
\end{table}

Examining the table, we see that there are three consecutive numbers,
6, 7, 8, which are necessarily followed by three consecutive odd numbers,
$c(6)=9$, $c(7)=11$, $c(8)=13$.
Thus 9 is present, 10 is missing, 11 is present, 12 is missing
and 13 is present.
Therefore the sequence continues with $c(9) =15$ (odd),
$c(10) = 16$ (even), $\ldots$, $c(13)=19$ (odd), $c(14) = 20$ (even).
This behavior is repeated forever.
A run of consecutive numbers is immediately followed by
a run of the same length of consecutive odd numbers.
And a run of consecutive odd numbers is immediately
followed by a run of twice that 
length of consecutive numbers (alternating even and odd).
Once we have noticed this, it is straightforward to find
an explicit formula that describes this sequence:
$$
c(1) = 1, \quad c(2) =4 ,
$$
and subsequent terms are given by
$$
c(9 \cdot 2^k - 3+j) = 12 \cdot 2^k -3 + \frac{3}{2} j + \frac{1}{2} |j|
$$
for $k \ge 0$, $-3\cdot 2^k \le j < 3\cdot 2^k$ (see \cite{CSV03} for the proof).

The structure is further revealed by
examining the sequence of first differences,
$c(n) = c(n+1) - c(n)$, which is
$$
3,2,1,1,1,2,2,2,1^6, 2^6,1^{12}, 2^{12}, 1^{24}, 2^{24}, \ldots \,,
$$
where we have written $1^6$ to indicate a run of six 1's, etc.
The oscillations double in length at each step.

\paragraph{4. Approximate squaring.}
The symbol
$\lceil x \rceil$ denotes the ``ceiling'' function,
the smallest integer greater than or equal to $x$.
Start with any fraction greater than $1$, say
$\F{8}{7}$, and repeatedly apply the ``approximate squaring'' map:
\beql{Eq2}
\mbox{replace~} x \mbox{~by~} x \lceil x \rceil \,.
\eeq
Since $\F{8}{7} = 1.142\ldots$, $\lceil \F{8}{7} \rceil = 2$,
so after the first step we reach $\F{8}{7} \times 2 = \F{16}{7}$.
A second approximate squaring step takes us to 
$\F{16}{7} \times 3 = \F{48}{7}$, 
and a third step takes us to 
$\F{48}{7} \times 7 = 48$,
which is an integer, and we stop. It took three steps to reach an integer.
The question is: {\em do we always reach an integer}?
Jeffrey Lagarias and I studied this problem in \cite{APSQ}.
We showed that almost all initial fractions greater than $1$
eventually reach an integer,
and that if the denominator is $2$ then they all do,
but we were unable to give an affirmative answer in general.
In fact, we show that the problem has some similarities to the 
notorious Collatz (or ``$3x+1$'')
problem \cite{Lagarias1985}, and so may be difficult
to solve in general.
(A similar problem has been posed by Jim Tanton \cite{Tanton2002}.)

The numbers involved grow very rapidly: if we start
with $\F{6}{5}$, for example, successive approximate squarings
produce the sequence
\begin{align}\label{Eq4}
& {\frac{6}{5}}\,,
 {\frac {12}{5}}\,,
 {\frac {36}{5}}\,,
 {\frac {288}{5}}\,,
 {\frac {16704}{5 }}\,,  
 {\frac {55808064}{5}}\,,
{\frac {622908012647232}{5}}\,, \nonumber \\
& {\frac { 77602878444025201997703040704}{5}}\,, \nonumber \\
& {\frac { 1204441348559630271252918141028336694332989128001036771264}{5}}\,,
\ldots \nonumber
\end{align}
(cf. A117596), which finally reaches an integer, a $57735$-digit number,
after $18$ steps!

If the fraction that we start with has denominator $2$, we can say
exactly how many steps are needed.
If we start with $\F{2l+1}{2}$ then we reach an integer
in $m+1$ steps, where $2^m$ is the highest power of $2$ that divides $l$.
For example, $\F{17}{2}$ (where $l=2^3$) reaches 
the integer $1204154941925628$ in $4$ steps.

But even for denominator $3$ we cannot say exactly what will happen.
The following table shows what happens for the first few starting
values.
It gives the initial term, the number of steps to reach an integer,
and the integer that is reached.
$$
\begin{array}{rcccccccccc}
\mbox{start}: & \frac{3}{3} & \frac{4}{3} & \frac{5}{3} & \frac{6}{3} & \frac{7}{3} & \frac{8}{3} & \frac{9}{3} & \frac{10}{3} & \frac{11}{3} & \cdots \\
\mbox{steps}: & 0 & 2 & 6 & 0 & 1 & 1 & 0 & 5 & 2 & \cdots \\
\mbox{reaches}: & 1 & 8 & 1484710602474311520 & 2 & 7 & 8 & 3 & 1484710602474311520 & 220 & \cdots
\end{array}
$$
(The second and third rows are sequences 
A072340 and A085276---they certainly stagger.)

Starting values of the form $\F{l+1}{l}$ seem to
take an especially long time to reach an integer.
The examples $l=5$ and $7$ have already been mentioned.
It is amusing to note that if we start with
$\frac{200}{199}$ and
repeatedly apply the approximately squaring operation,
the first integer reached is roughly
$$
200^{2^{1444}} ~,
$$
a number with about $10^{435}$ digits.

\paragraph{5. ``If a power series was a power of a power series,
what power would it be, seriously?''}
That was the title of Nadia Heninger's talk about her 2005  summer project
at AT\&T Labs.
Consider the following question, a typical problem from
number theory.
How many ways are there to write a given number $n$ as the 
sum of four squares?
That is, how many integer solutions $(i,j,k,l)$ are there
to the equation
$$
n ~=~ i^2+j^2+k^2+l^2 ~?
$$
Call the answer $r_4(n)$.
Solutions with $i,j,k,l$ in a different order or
with different signs count as different,
so for instance $r_4(4) = 24$,
since for $n=4$ $(i,j,k,l)$ can be any of
$$
(\pm2,0,0,0), (0,\pm2,0,0), (0,0,\pm2,0), (0,0,0,\pm2),
(\pm1,\pm1,\pm1,\pm1) \,.
$$
We can capture this problem in a generating 
function that looks like
$$
R(x) := 
r_4(0) + r_4(1)x + r_4(2)x^2 + r_4(3)x^3 + \cdots \,,
$$
in which the coefficient of $x^n$ gives the answer $r_4(n)$.
For this problem, it is easy to see that $R(x)$
is equal to the fourth power of Jacobi's famous 
``theta series''
$$
\theta_3(x) :=  1 + 2x+ 2x^4 + 2x^9 + 2x^{16} + 2x^{25} + \cdots   
$$
(this is classical number theory: see for example 
Hardy and Wright \cite{HW} or Grosswald \cite{Gross}).
Of course this means that we can take the fourth root of 
$R(x)$ and still have integer coefficients.

What we (that is, Nadia Heninger, Eric Rains and I \cite{OSCE})
discovered is that there are many other
important generating functions for which it is possible
to take a fourth root, or in general a $k^{\rm th}$ root,
and still have integer coefficients.

Many of our examples arise as ``theta functions'' of sphere packings.
The most familiar sphere packing is the grocer's
face-centered cubic lattice arrangement of oranges,
which Tom Hales \cite{Hales} has recently shown to be the densest possible
sphere-packing of three-dimensional balls.
For any lattice packing of balls in $N$-dimensional space,
the theta series is a generating function 
whose coefficients give the number of balls with centers
at a given distance from the origin.
If there are  $M_d$ balls whose distance from the origin
is $\sqrt{d}$, the theta series is
\beql{EqTH}
\sum_{d} M_d x^d \,.
\eeq
The Jacobi theta series $\theta_3(x)$ mentioned above
is simply the theta series of the one-dimensional
lattice formed by the integers, and $R(x)$ is the theta series of the 
simple cubic lattice in four dimensions.

The starting point for our work was an observation of
Michael Somos \cite{Somos05} that
the $12$-th root of the theta series of
a certain $24$-dimensional lattice discovered by
Gabriele Nebe (\cite{Nebe98}, sequence A004046),
appeared to have integer coefficients.
We were able to establish his conjecture,
and to generalize it to many other theta series and power series.

An example of one of our discoveries is this: the
theta series of the densest lattice sphere packing in four
dimensions is the fourth power of a generating
function with integer coefficients.
The theta series in question, that of the $D_4$ lattice packing, is
\beql{EqD4}
r_4(0) + r_4(2)x^2 + r_4(4)x^4 + \cdots 
~=~ 1 + 24x^2 + 24x^4+96x^6+24x^8+144x^{10}+\cdots \,,
\eeq
and is formed by taking the even powers of $x$ in $R(x)$.
When we take the fourth root, we get
\begin{eqnarray*}
& 1 & + ~ 6\,{x}^{2}-48\,{x}^{4}+672\,{x}^{6}-10686\,{x}^{8}+185472\,{x}^{10}-3398304\,{x}^{12} \\
&&{} ~ +64606080\,{x}^{ 14} -1261584768\,{x}^{16}+25141699590\,{x}^{18}-\cdots \,.
\end{eqnarray*}
The coefficients stagger, changing sign
at each step, and growing in size (A108092).
Can a reader find any other interpretation 
of these coefficients?

In our paper, we begin by
studying the more general question of when a power
series of the form 
$$
F(x) = 1 + f_1x+f_2x^2+f_3x^3+f_4x^4+\cdots 
$$
with integer coefficients is the $k^{\rm th}$ power
of another such power series.
One of our results is that $F(x)$
is a $k^{\rm th}$ power if and only if
the series obtained by reducing the coefficients
of $F(x)$ mod $\mu_k$ is a $k^{\rm th}$ power,
where $\mu_k$ is obtained by multiplying $k$ by 
all the distinct primes that divide it.
E.g., to test if $F(x)$ is a fourth power, it is enough to
check the series obtained by
reducing the coefficients mod $\mu_4 = 4 \times 2 = 8$.
Since all the coefficients of \eqn{EqD4}
except the constant term are divisible by $8$,
it reduces to $1$ mod $8$, and certainly $1$ {\em is}
a fourth power of a series with integer coefficients! So \eqn{EqD4} is also.

More dramatic examples are provided by the theta series 
of the $E_8$ lattice in $8$ dimensions
and the Leech lattice in $24$ dimensions,
which are respectively $8^{\rm th}$ and
$24^{\rm th}$ powers of series with integral (albeit staggering)
coefficients.

In \cite{OSCE} we give many other examples,
including generating functions arising from weight enumerators of codes.

\paragraph{6. Dissections.}
My colleague Vinay Vaishampayan and I have found some surprising
applications of the classical dissection problem to optical communications.
But it seems that even the simplest questions in this subject are 
still unanswered.  Since
many Martin Gardner fans are experts in puzzles and their history,
I mention this here, in the hope
that some reader can provide further information.

\begin{figure}[htb]
\centerline{\includegraphics[ scale=.50 ]{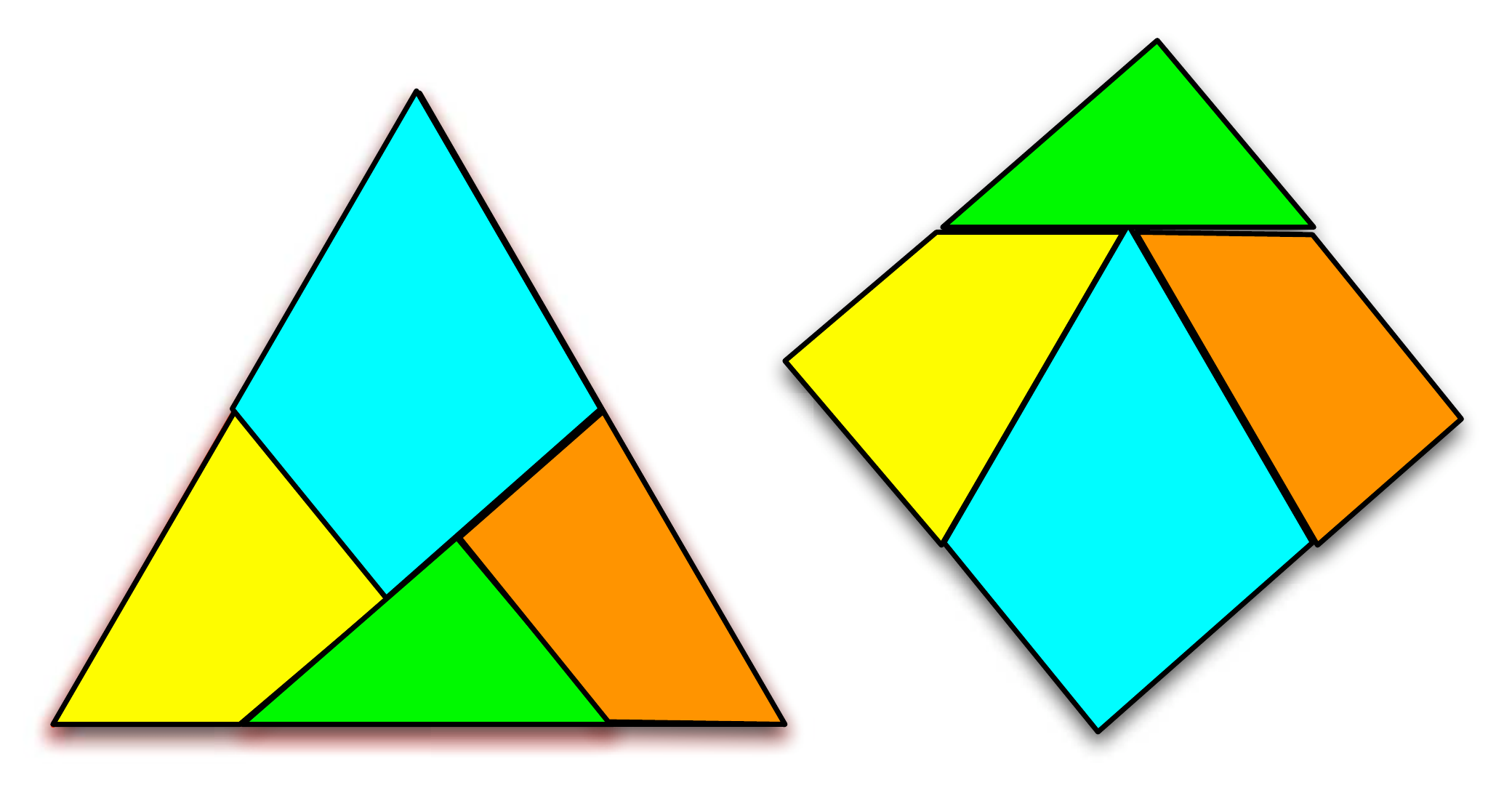}}

\caption{A triangle can be cut into four pieces
which can be rearranged to form a square.
Is it known that this cannot be done using only three pieces?}
\label{FVV}
\end{figure}

The following is the
simplest version of the question.
It is known that any polygon can be cut up
into a finite number of pieces
which can be arranged, without overlapping, to
form a square of the same area. (You are allowed to
turn pieces over. To avoid complications
caused by non-measurable sets, the edges of the pieces
must be simple curves.)
The question is, what is the minimal number, $d(n)$ say,
of pieces that are required to dissect a regular
polygon with $n$ sides ($n \ge 3$) into a square?
For the case $n=3$, we are looking for a minimal dissection of 
an equilateral triangle to a square. 
There is a famous $4$-piece dissection, apparently first
published by Dudeney in 1902 
(\cite{Dudeney02},  \cite{Fred}), 
shown in Fig. \ref{FVV}.
It seems unlikely that a three-piece dissection exists,
but has anyone ever proved this?
In other words, is $d(3)$ really $4$?

As far as I know, none of the values of $d(n)$ are known
for certain (except of course $d(4)=1$).
The best values presently known for $d(3), d(4), d(5), \ldots$
(A110312), taken from Frederickson \cite{Fred}
and Theobald \cite{Theo}, are:
$$
4?, 1, 6?, 5?, 7?, 5?, 9?, 7?, \ldots \,.
$$
This is most unsatisfactory: the normal rule is
that every term in a sequence in the OEIS should be
known to be correct. This sequence is quite an exception,
the values shown being merely upper bounds.
It is not surprising that these entries stagger a bit---but
are they correct?

\begin{figure}[htb]
\centerline{\includegraphics[ width=2in]{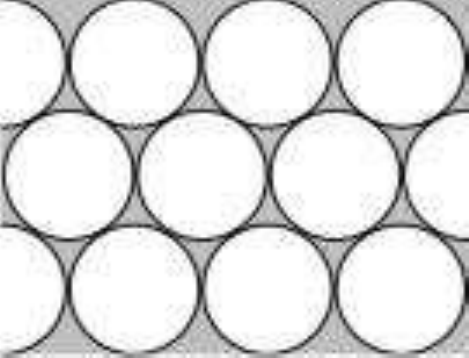}}
\caption{A portion of 
the hexagonal lattice packing of circles in the plane,
illustrating the solution to
the two-dimensional kissing number problem.}
\label{FigPenny}
\end{figure}

\paragraph{7. The kissing number problem.}
The $N$-dimensional kissing number problem 
asks for the maximal number of
$N$-dimensional balls that can touch another ball of the same radius
(the term comes from billiards).
The problem has applications in geometry, number theory, group theory
and digital communications (\cite{SPLAG}, \cite{ERS}).
For example, in two dimensions, six pennies
is the maximal number than can touch another penny, as shown in
Fig. \ref{FigPenny}. This illustration shows a portion of
the familiar hexagonal lattice packing in the plane.
For the solution of the problem, however, 
the balls need not necessarily be part of a lattice packing.
In dimensions $1$ to $4$, $8$ and $24$ the highest possible kissing number
can be achieved  using a lattice packing,
but in dimension $9$ there is
a nonlattice packing with a maximal kissing number that is higher
than is possible in any lattice packing,
and this is almost certainly also true
in ten dimensions---see \cite{SPLAG}.
If the record is achieved by
a lattice packing, then we can read off the kissing
number from the theta series defined in \eqn{EqTH}: 
this begins  $1 + \tau x^u + \cdots$
for some $u$, where $\tau$ is the kissing number.

The answer in three dimensions is $12$, and Musin \cite{Mus}, \cite{Mus2}
has recently established the long-standing conjecture
that the answer in four dimensions is $24$,
as found in the $D_4$ lattice packing---the number can
be read off the theta series in \eqn{EqD4}.  We also
know the answers in dimensions eight and twenty-four
($240$ and $196560$, respectively).
The beginning of the sequence
of solutions to the $N$-dimensional kissing number problem is
$$
2, 6, 12, 24, 40?, 72?, 126?, 240, 306?, 500?, \ldots , 196560, \ldots 
$$
(cf. A001116). The entries with question marks are
merely lower bounds.
If an oracle offered to supply $64$ terms of any sequence that I chose,
I would pick this one.
I would also ask the oracle for the constructions that it used---in particular,
in high dimensions, are the arrangements always without structure,
or is there an infinite sequence of dimensions where there are 
elegant algebraic constructions?

\end{document}